\documentclass{amsart}
\usepackage{amssymb}
\usepackage{amsmath}
\usepackage{amsfonts}

\newtheorem{theorem}{Theorem}
\theoremstyle{plain}
\newtheorem{acknowledgement}{Acknowledgement}

\newtheorem{condition}{Condition}

\newtheorem{lemma}{Lemma}

\newtheorem{remark}{Remark}

\numberwithin{equation}{section}
\input{tcilatex}

\begin{document}
\title{Smooth densities for Stochastic Differential Equations with Jumps}
\author{THOMAS CASS }
\date{October 2, 2007}
\maketitle

\begin{abstract}
We consider a solution $x_{t}$ to a generic Markovian jump diffusion and
show that for any $t_{0}>0$ the law of $x_{t_{0}}$ has a $C^{\infty }$
density with respect to Lebesgue measure under a uniform version of H\"{o}%
rmander's conditions. Unlike previous results in the area the result covers
a class of infinite activity jump processes. The result is accompolished by
using carefully crafted refinements to the classical arguments used in
proving the smoothness of density via Malliavin calculus. In particular, we
provide a proof that the semimartinagale inequality of Norris persists for
discontinuous semimartingales when the jumps are small.
\end{abstract}

\section{Introduction}

This paper focuses on the study of the stochastic differential equation

\begin{equation}
x_{t}=x+\int_{0}^{t}Z(x_{s-})ds+\int_{0}^{t}V(x_{s-})dW_{s}+\int_{0}^{t}%
\int_{E}Y(x_{s-},y)(\mu -\nu )(dy,ds),  \label{SDE}
\end{equation}

\noindent and addresses the fundamental problem of finding a sufficient
condition for the existence of a smooth ($C^{\infty }$) density for the
solution at positive times. For diffusion processes the pioneering work of
Bismut \cite{bismut} \ and Stroock \cite{stroock1} and \cite{stroock2}
provides a probabilistic framework for establishing such a result under H%
\"{o}rmander's conditions on the vector fields. As is pointed out in \cite%
{stroock2} it is, given the existence of alternative methods based on
partial differential equations, difficult to justify the effort involved in
the probabilistic proof of this result purely for the sake of diffusion
processes. From the outset it was always understood that this approach
should be used as a template for investigating the smoothness properties for
different probabilistic objects, not amenable to analysis by PDE theory. We
now switch our focus to the question : \textit{when does a solution to the
SDE (\ref{SDE}) admit a smooth density?}

We point out that we are by no means the first to consider this problem and\
several prominent landmarks are worthy of comment. \ The first comprehensive
account of these ideas was presented in \cite{bichteler2}, where a
smoothness result is proved under a uniform ellipticity on the diffusion
vector fields (in fact \cite{bichteler2} also explores how a smooth density
can be acquired through the jump component)$.$ Further progress was made in 
\cite{norris1} where existence of the density was shown under a version of H%
\"{o}rmander's conditions which are local in the starting point. Both these
works were successful in establishing a criterion for a smooth density
namely that the inverse of \ the (reduced) Malliavin covariance matrix has
finite $L^{p}$ norms for $p\geq 2$.

Verification of this criterion usually occurs by way of subtle estimates on
the reduced covariance matrix which are in general difficult to establish.
In the diffusion case a streamlined approach to obtaining these estimates
has been achieved by a semimartingale inequality known as Norris's lemma
(see \cite{norris} or \cite{Nu06}). This result, interesting in its own
right, provides an estimate for the probability that a continuous
semimartingale is small on a set where its quadratic variation is
comparatively large. Traditionally, this result has been presented as a
quantitative form of the uniqueness of the Doob-Meyer decomposition for
continuous semimartingales, however the appearance of similar estimates in
the context of fractional Brownian motion with $H>1/2$ (not a
semimartingale) (see \cite{BH}) have made it seem as though Norris's lemma
expresses something fundamental rather than anything tied to the particular
structure of continuous semimartingales.

Some recent work in the case of jump diffusions has been undertaken in \cite%
{teichmann},\cite{komatsu} and \cite{tak}. The article \cite{teichmann}
proves a smoothness result under uniform H\"{o}rmander conditions and under
the assumption that the underlying jump process is of finite activity. This
is achieved by fixing some $T>0$, conditioning on $N_{T}=n$, the number of
jumps until time $T,$ and noticing that this gives rise to some (random)
interval $[S_{1}(\omega ),S_{2}(\omega ))$ with $0\leq S_{1}<S_{2}<T$ such
that $S_{2}(\omega )-S_{1}(\omega )\geq T(n+1)^{-1}$ and 
\begin{equation*}
\left\{ x_{t}^{x}:S_{1}\leq t<S_{2}\right\} \overset{\mathcal{D}}{=}\left\{ 
\tilde{x}_{t}^{x_{S_{1}}^{x}}:0\leq t<S_{2}-S_{1}\right\}
\end{equation*}

\noindent where\bigskip\ $\tilde{x}_{t}^{x}$ is the diffusion process%
\begin{equation*}
\tilde{x}_{t}^{x}=x+\int_{0}^{t}Z(\tilde{x}_{s}^{x})ds+\int_{0}^{t}V(\tilde{x%
}_{s}^{x})dW_{s}.
\end{equation*}%
\noindent The usual diffusion Norris lemma may be applied to give estimates
for the Malliavin covariance matrix arising from $\tilde{x}_{t}$ on this
interval which can then easily be related to covariance matrix for $x_{t}$.
In this paper we pursue this idea further by proving that the quality of the
estimate which features in Norris's lemma is preserved when jumps are
introduced provided that these jumps are small enough that they do not
interfere too much. \ We then develop the conditioning argument outlined
above by splitting up the sample path into disjoint intervals on which the
jumps are small, and then estimating the Malliavin covariance matrix on the
largest of these intervals. The outcome of this reasoning will be the
conclusion that a solution to (\ref{SDE}) has a smooth density under uniform
H\"{o}rmander conditions (indeed, the same conditions as in \cite{teichmann}%
) and subject to some restrictions on the rate at which the jump measure
accumulates small jumps. These conditions are sufficiently flexible to admit
some jump diffusions based on infinite activity jump processes.

This paper is arranged as follows : we first present some preliminary
results and notation on Malliavin calculus. Subsequently, we state and prove
our new version of Norris's lemma and then illustrate how it may be utilized
in concert with classical arguments to verify the $C^{\infty }$ density
criterion for the solution to (\ref{SDE}).

\begin{acknowledgement}
The author would like to thank James Norris and Peter Friz for related
discussions.
\end{acknowledgement}

\section{\protect\bigskip Preliminaries}

\bigskip Let $x_{t}$ denote the solution to the SDE%
\begin{equation}
x_{t}=x+\int_{0}^{t}Z(x_{s-})ds+\int_{0}^{t}V(x_{s-})dW_{s}+\int_{0}^{t}%
\int_{E}Y(x_{s-},y)(\mu -\nu )(dy,ds),  \label{sol}
\end{equation}

\noindent Where $W_{t}=\left( W_{t}^{1},...,W_{t}^{d}\right) $ is an $%
%TCIMACRO{\U{211d} }%
%BeginExpansion
\mathbb{R}
%EndExpansion
^{d}-$valued Brownian motion on some probability space $\left( \Omega ,%
\mathcal{F}_{t},P\right) $ and $\mu $ is a $\left( \Omega ,\mathcal{F}%
_{t},P\right) -$Poisson random measure on $E\times \lbrack 0,\infty )$ for
some topological \footnote{%
we will later need some vector space structure on $E$ and will principally
be concerned with the case $E=\mathbb{R}^{n}$} space $E$ such that $\nu $,
the compensator of $\mu $, is of the form $G(dy)dt$ for some $\sigma $%
-finite measure $G.$ The vector fields $Z:%
%TCIMACRO{\U{211d} }%
%BeginExpansion
\mathbb{R}
%EndExpansion
^{e}\rightarrow 
%TCIMACRO{\U{211d} }%
%BeginExpansion
\mathbb{R}
%EndExpansion
^{e},Y\left( \cdot ,y\right) :%
%TCIMACRO{\U{211d} }%
%BeginExpansion
\mathbb{R}
%EndExpansion
^{e}\rightarrow 
%TCIMACRO{\U{211d} }%
%BeginExpansion
\mathbb{R}
%EndExpansion
^{e}$ and $V=(V_{1},...,V_{d}),$ where $V_{i}:%
%TCIMACRO{\U{211d} }%
%BeginExpansion
\mathbb{R}
%EndExpansion
^{e}\rightarrow 
%TCIMACRO{\U{211d} }%
%BeginExpansion
\mathbb{R}
%EndExpansion
^{e}$ for $i\in \left\{ 1,...,d\right\} .$ At times we will write $x_{t}^{x}$
when we wish to emphasize the dependence of the process on its initial
condition.

We introduce some notation, firstly for $p\in 
%TCIMACRO{\U{211d} }%
%BeginExpansion
\mathbb{R}
%EndExpansion
$ let 
\begin{equation*}
L_{+}^{p}(G)=\left\{ f:E\rightarrow 
%TCIMACRO{\U{211d} }%
%BeginExpansion
\mathbb{R}
%EndExpansion
^{+}:\int_{E}f(y)^{p}G(dy)<\infty \right\} ,
\end{equation*}

\noindent and define%
\begin{equation*}
L_{+}^{p,\infty }(G)=\underset{q\geq p}{\cap }L_{+}^{q}(G).
\end{equation*}

\noindent \noindent We will always assume that at least the following
conditions are in force

\begin{condition}
\label{c1}$Z,V_{1},...,V_{d}\in C_{b}^{\infty }($ $%
%TCIMACRO{\U{211d} }%
%BeginExpansion
\mathbb{R}
%EndExpansion
^{e})$
\end{condition}

\begin{condition}
\label{c2}For some $\rho _{2}\in L_{+}^{2,\infty }(G)$ and every $n\in 
%TCIMACRO{\U{2115} }%
%BeginExpansion
\mathbb{N}
%EndExpansion
$%
\begin{equation*}
\underset{y\in E,x\in 
%TCIMACRO{\U{211d} }%
%BeginExpansion
\mathbb{R}
%EndExpansion
^{e}}{\sup }\frac{1}{\rho _{2}(y)}|D _{1}^{n}Y(x,y)|<\infty .
\end{equation*}
\end{condition}

\begin{condition}
\label{c3}$\sup_{y\in E,x\in 
%TCIMACRO{\U{211d} }%
%BeginExpansion
\mathbb{R}
%EndExpansion
^{e}}|\left( I+D _{1}Y(x,y)\right) ^{-1}|<\infty .$
\end{condition}

\noindent We now define the processes $J_{t\leftarrow 0}$ and $%
J_{0\leftarrow t}$ considered as linear maps from $%
%TCIMACRO{\U{211d} }%
%BeginExpansion
\mathbb{R}
%EndExpansion
^{e}$ to $%
%TCIMACRO{\U{211d} }%
%BeginExpansion
\mathbb{R}
%EndExpansion
^{e}$ as the solutions to the following SDEs 
\begin{eqnarray}
J_{t\leftarrow 0} &=&I+\int_{0}^{t}DZ(x_{s-})J_{s-\leftarrow
0}ds+\int_{0}^{t}DV(x_{s-})J_{s-\leftarrow 0}dW_{s}  \label{jacobian} \\
&&+\int_{0}^{t}\int_{E}D_{1}Y(x_{s-},y)J_{s-\leftarrow 0}(\mu -\nu )(dy,ds) 
\notag
\end{eqnarray}

\noindent and%
\begin{gather}
J_{0\leftarrow t}=I-\int_{0}^{t}J_{0\leftarrow s-}\Bigg(DZ(x_{s-})-\underset{%
i=1}{\overset{d}{\sum }}DV_{i}(x_{s-})^{2}  \label{ijacobian1} \\
-\int_{E}(I+D_{1}Y(x_{s-},y))^{-1}D_{1}Y(x_{s-},y)^{2}G(dy)\Bigg)%
ds-\int_{0}^{t}J_{0\leftarrow s-}DV(x_{s-})dW_{s}  \notag \\
-\int_{0}^{t}\int_{E}J_{0\leftarrow
s-}(I+D_{1}Y(x_{s-},y))^{-1}D_{1}Y(x_{s-},y)(\mu -\nu )(dy,ds).  \notag
\end{gather}

\noindent The following result may then be verified (see for instance \cite%
{norris})

\begin{theorem}
\label{jacob}Under conditions \ref{c1},\ref{c2} and \ref{c3} the system of
SDEs (\ref{sol},\ \ref{jacobian}) and \ (\ref{sol},\ref{ijacobian1}) have
unique solutions with 
\begin{equation*}
\underset{0\leq s\leq t}{\sup }|J_{s\leftarrow 0}|\text{ and}\underset{0\leq
s\leq t}{\sup }|J_{0\leftarrow s}|\text{ }\in L^{p}
\end{equation*}

\noindent for all $t\geq 0$ and $p<\infty .$ Moreover,%
\begin{equation*}
J_{0\leftarrow t}=J_{t\leftarrow 0}^{-1}\text{ for all }t\geq 0\text{ almost
surely.}
\end{equation*}
\end{theorem}

\noindent We define the reduced Malliavin covariance matrix%
\begin{equation*}
C_{0,t}^{x,I}=C_{t}^{x,I}=\int_{0}^{t}\underset{i=1}{\overset{d}{\sum }}%
J_{0\leftarrow s-}^{x,I}V_{i}(x_{s-}^{x})\otimes J_{0\leftarrow
s-}^{x,I}V_{i}(x_{s-}^{x})ds
\end{equation*}

\noindent which we will sometimes refer to simply as $C_{t}$ suppressing the
dependence on the initial conditions. The following well known result
provides a sufficient condition for the process $x_{t}$ to have a $C^{\infty
}$ density in terms of the moments of the inverse of $C_{t}.$

\begin{theorem}
\label{criterion}Fix $t_{0}>0$ and $x\in 
%TCIMACRO{\U{211d} }%
%BeginExpansion
\mathbb{R}
%EndExpansion
^{e}$ and suppose that for every $p\geq 2$ $\left\vert C_{t}^{-1}\right\vert
\in L^{p},$ then $x_{t_{0}}^{x}$ has a $C^{\infty }$density with respect to
Lebesgue measure.
\end{theorem}

\section{Norris's lemma}

From now on we set $E=\mathbb{R}^{n}$. The following result provides an
exponential martingale type inequality for a class of local martingales
based on stochastic integrals with respect to a Poisson random measure when
the jumps of the local martingale are bounded. Interesting discussions on
results of this type can be found in \cite{barlow} and \cite{dzh}

\begin{lemma}
\label{EMI}\noindent Let $\mu $ be a Poisson random measure on $E\times
\lbrack 0,\infty )$ with compensator $\nu $ of the form $\nu
(dy,dt)=G(dy)dt. $ Let $f(t,y)$ be a real-valued previsible process having
the property that 
\begin{equation*}
\underset{y\in E}{\sup }\text{ }\underset{0\leq s\leq t}{\sup }|f(s,y)|<A
\end{equation*}

\noindent for every $0<t<\infty $ and some $A<\infty $. Then, if $%
M_{t}=\int_{0}^{t}\int_{E}f(s,y)(\mu -\nu )(dy,ds)$ the following inequality
holds%
\begin{equation*}
P\left( \underset{0\leq s\leq t}{\sup }|M_{s}|\geq \delta ,\left\langle
M\right\rangle _{t}<\rho \right) \leq 2\exp \left( -\frac{\delta ^{2}}{%
2(A\delta +\rho )}\right)
\end{equation*}
\end{lemma}

\begin{proof}
Consider $Z_{t}=\exp (\theta M_{t}-\alpha \left\langle M\right\rangle _{t})$
with $0<\theta <A^{-1}$ and $\alpha =2^{-1}\theta ^{2}(1-\theta A)^{-1}.$
Since for any $x\in 
%TCIMACRO{\U{211d} }%
%BeginExpansion
\mathbb{R}
%EndExpansion
$ we have 
\begin{equation}
g_{\theta }(x):=e^{\theta x}-1-\theta x=\underset{k=2}{\overset{\infty }{%
\sum }}\frac{\theta ^{k}x^{k}}{k!}\leq \frac{\theta ^{2}x^{2}}{2}\underset{%
k=0}{\overset{\infty }{\sum }}(\theta A)^{k}=\frac{\theta ^{2}x^{2}}{%
2(1-\theta A)}=\alpha x^{2}.  \label{expmart}
\end{equation}

\noindent We may deduce that $Z$ is a supermartingale by writing 
\begin{eqnarray*}
Z_{t} &=&\exp \left( \theta M_{t}-\int_{0}^{t}\int_{E}g_{\theta
}(f(s,y))G(dy)ds\right) \\
&&\exp \left( \int_{0}^{t}\int_{E}\left( g_{\theta }(f(s,y))-\alpha
f(s,y)^{2}\right) G(dt)ds\right)
\end{eqnarray*}

\noindent and, using It\^{o}'s formula the first term of the product is a
non-negative local martingale (and hence a supermartingale) and the second
term decreases in $t$ by (\ref{expmart}). Define the stopping time $T=\inf
\left\{ s\geq 0:\left\langle M\right\rangle _{s}>\rho \right\} $ then, since 
$E[Z_{T}]\leq 1,$ taking $\theta =\delta (\rho +A\delta )^{-1}$ and applying
Chebyshev's inequality gives%
\begin{equation*}
P\left( \underset{0\leq s\leq t}{\sup }|M_{s}|\geq \delta ,\left\langle
M\right\rangle _{t}<\rho \right) \leq P\left( \underset{0\leq s\leq T}{\sup }%
Z_{s}\geq e^{\delta \theta -\alpha \rho }\right) \leq \exp \left( -\frac{%
\delta ^{2}}{2(A\delta +\rho )}\right) .
\end{equation*}

\noindent Finally, we complete the proof by applying the same argument to $%
-M.$
\end{proof}

\bigskip

\noindent From now on we will assume that the following technical conditions
on the jump measure $G$ and the jump vector field $Y$ are in force :

\begin{condition}
\label{condition1}$\underset{x\in 
%TCIMACRO{\U{211d} }%
%BeginExpansion
\mathbb{R}
%EndExpansion
^{e}}{\sup }\int_{E}|Y(x,y)|G(dy)<\infty .$
\end{condition}

\begin{condition}
\label{condition2}For some $\kappa \geq n$ we have 
\begin{equation}
\underset{\epsilon \downarrow 0}{\lim \sup }\frac{1}{f(\epsilon )}%
\int_{|y|>\epsilon }G(dy)\text{ }<\infty \text{ ,}  \label{cond2a}
\end{equation}

\noindent where $f:(0,\infty )\rightarrow (0,\infty )$ is defined by 
\begin{equation}
f(x)=\left\{ 
\begin{array}{c}
\log x^{-1}\text{\ if }\kappa =n \\ 
x^{-\kappa +n}\text{\ \ if }\kappa >n%
\end{array}%
\right. .  \label{cond2b}
\end{equation}

\noindent Moreover, for any $\beta >0$ we have%
\begin{equation*}
\int_{E}|y|^{\kappa -n+\beta }G(dy)<\infty ,
\end{equation*}

\noindent and%
\begin{equation}
\underset{\epsilon \downarrow 0}{\lim \sup }\frac{1}{\epsilon ^{\beta }}%
\int_{|y|<\epsilon }|y|^{\kappa -n+\beta }G(dy)\text{ }<\infty \text{ .}
\label{cond2c}
\end{equation}
\end{condition}

\begin{condition}
\label{condition3}There exists a function $\phi \in L_{+}^{1}(G)$ which has
the properties that for some $\alpha >0$ 
\begin{equation*}
\underset{y\rightarrow 0}{\lim \sup }\frac{\phi (y)}{|y|^{\kappa -n+\alpha }}%
\text{ }<\infty ,
\end{equation*}

\noindent and, for some positive constant $C<\infty $ and every $k\in 
%TCIMACRO{\U{2115} }%
%BeginExpansion
\mathbb{N}
%EndExpansion
\cup \left\{ 0\right\} $%
\begin{equation*}
\underset{x\in 
%TCIMACRO{\U{211d} }%
%BeginExpansion
\mathbb{R}
%EndExpansion
^{e}}{\sup }|D _{1}^{k}Y(x,y)|\leq C\phi (y).
\end{equation*}
\end{condition}

\bigskip

Conditions \ref{condition1}, \ref{condition2} and \ref{condition3} may at
first sight appear somewhat opaque, however they will be crucial ingredient
in our subsequent arguments, in particular they enable us quantify the rate
at which the total mass of the jump measure increases near zero. \ To
develop intuition for their implications consider the following
straight-forward example : take $n=1$ and $Y(x,y)=\widetilde{Y}(x)y$ for
some $C^{\infty }-$bounded $\widetilde{Y}:%
%TCIMACRO{\U{211d} }%
%BeginExpansion
\mathbb{R}
%EndExpansion
^{e}\rightarrow 
%TCIMACRO{\U{211d} }%
%BeginExpansion
\mathbb{R}
%EndExpansion
^{e}$ (this puts us in the set up of \cite{teichmann})$.$ Also, define the
measure $G$ on $%
%TCIMACRO{\U{211d} }%
%BeginExpansion
\mathbb{R}
%EndExpansion
$ by taking $G(dy)=|y|^{-\kappa }dy$. We then see what is needed to verify
each of the conditions in turn, firstly, condition \ref{condition1} will be
satisfied provided 
\begin{equation*}
\underset{x\in 
%TCIMACRO{\U{211d} }%
%BeginExpansion
\mathbb{R}
%EndExpansion
^{e}}{\sup }\int_{E}|\widetilde{Y}(x)y|G(dy)=\underset{x\in \mathbb{R}^{e}}{%
\sup}|\widetilde{Y}(x)| \int_{E}|y|G(dy)=2\underset{x\in \mathbb{R}^{e}}{\sup%
}|\widetilde{Y}(x)| \int_{0}^{\infty}y^{1-\kappa }dy<\infty ,
\end{equation*}

\bigskip \noindent which will hold so long as $\kappa <2$. The constraint
that $\kappa \geq 1$ in condition \ref{condition2} ensures that the jump
measure is of infinite activity and (\ref{cond2a}) and (\ref{cond2b}) are
trivially verified by integration. Since we are in the setting $1\leq \kappa
<2$, we may find $\alpha \in (0,1)$ such that $\kappa +\alpha <2$ to ensure
that $\phi (y):=|y|$ is $O(|y|^{\kappa -n+\alpha })$ as $y\rightarrow 0$ and
hence condition \ref{condition3} is also satisfied$.$

Suppose now that $\Upsilon :[0,t_{0}]\times E\rightarrow 
%TCIMACRO{\U{211d} }%
%BeginExpansion
\mathbb{R}
%EndExpansion
.$ is some given, real-valued, previsible process. It will at times be
important for us to impose the following condition on $\Upsilon$ .

\begin{condition}
\label{condition4}Let $G$ satisfy condition \ref{condition2}. Then there
exists some previsible process $D_{t}$ taking values in $[0,\infty )$ with $%
\sup_{0\leq t\leq t_{0}}D_{t}\in L^{p}$ for all $p\geq 1,$ and a function $%
\phi $ $\in L_{+}^{1}(G)$ such that 
\begin{equation}
|\Upsilon (t,y)|\leq D_{t}\phi (y)\text{ \ \ for all }t\in \lbrack 0,t_{0}]%
\text{ and }y\in E,
\end{equation}

\noindent and for some $\alpha =\alpha _{\Upsilon }>0$ 
\begin{equation}
K_{\phi }:=\underset{y\rightarrow 0}{\lim \sup }\frac{\phi (y)}{|y|^{\kappa
-n+\alpha }}\text{ }<\infty .  \label{vfconstant}
\end{equation}
\end{condition}

Equipped with these remarks we are now in a position to state and prove the
following lemma which will be fundamental to providing the estimates on the
reduced covariance matrix we need later.

\begin{lemma}
\label{Norris}(Norris-type lemma) Fix $t_{0}>0$ and for every $\epsilon >0$
suppose $\beta^{\epsilon} (t),\gamma^{\epsilon} (t)=(\gamma^{\epsilon}
_{1}(t),...,\gamma^{\epsilon} _{d}(t))$ ,$u^{\epsilon}(t)=(u^{%
\epsilon}_{1}(t),...,u^{\epsilon}_{d}(t))$ are previsible processes taking
values in $%
%TCIMACRO{\U{211d} }%
%BeginExpansion
\mathbb{R}
%EndExpansion
,%
%TCIMACRO{\U{211d} }%
%BeginExpansion
\mathbb{R}
%EndExpansion
^{d}$ and $%
%TCIMACRO{\U{211d} }%
%BeginExpansion
\mathbb{R}
%EndExpansion
^{d}$ respectively. Suppose further that $\zeta^{\epsilon} (t,y)$and $%
f^{\epsilon}(t,y)$ are real-valued previsible processes satisfying condition %
\ref{condition4} such that the functions $\phi^{\zeta}$ and $\phi^{f}$ do
not depend on $\epsilon$ and moreover for every $q\geq 1$ 
\begin{equation}
\underset{\epsilon >0}{\sup }E\left[ \underset{0\leq t\leq t_{0}}{\sup }%
\left( D_{t}^{\zeta,\epsilon}\right) ^{q}+ \underset{0\leq t\leq t_{0}}{\sup 
}\left( D_{t}^{f,\epsilon}\right) ^{q}\right] <\infty.  \label{epsilonunif}
\end{equation}
Let $\alpha =\min (\alpha _{\zeta }$,$\alpha _{f})$, $\delta >0$, $z=3\delta
(\kappa -n+\alpha )^{-1}$ and define the processes $a^{\epsilon }$ and $%
Y^{\epsilon }$ as the solutions to the SDEs%
\begin{eqnarray*}
a^{\epsilon }(t) &=&\alpha +\int_{0}^{t}\beta^{\epsilon} (s)ds+\underset{i=1}%
{\overset{d}{\sum }}\int_{0}^{t}\gamma^{\epsilon}
_{i}(s)dW_{s}^{i}+\int_{0}^{t}\int_{|y|<\epsilon ^{z}}\zeta^{\epsilon}
(s,y)(\mu -\nu )(ds,dy) \\
Y^{\epsilon }(t) &=&y+\int_{0}^{t}a^{\epsilon}(s)ds+\underset{i=1}{\overset{d%
}{\sum }}\int_{0}^{t}u^{\epsilon}_{i}(s)dW_{s}^{i}+\int_{0}^{t}\int_{|y|<%
\epsilon ^{z}}f^{\epsilon}(s,y)(\mu -\nu )(ds,dy).
\end{eqnarray*}

\noindent Assume that for some $p\geq 2$ the quantity 
\begin{equation}
\underset{\epsilon >0}{\sup }E\left[ \underset{0\leq t\leq t_{0}}{\sup }%
\left( |\beta^{\epsilon} (t)|+|\gamma^{\epsilon} (t)|+|a^{\epsilon
}(t)|+|u^{\epsilon}(t)|+\int_{E}(|\zeta^{\epsilon}
(t,y)|^{2}+|f^{\epsilon}(t,y)|^{2})G(dy)\right) ^{p}\right]  \label{unif}
\end{equation}

\noindent is finite, and for some $\rho _{1},\rho _{2}\in L_{+}^{2,\infty
}(G)$ we have 
\begin{equation*}
\underset{\epsilon >0}{\sup } \left( E\left[ \left( \underset{0\leq t\leq
t_{0}}{\sup }\underset{y\in E}{\sup }\frac{|\zeta^{\epsilon}(t,y)|}{\rho
_{1}(y)}\right) ^{p}\right] +E\left[ \left( \underset{0\leq t\leq t_{0}}{%
\sup }\underset{y\in E}{\sup }\frac{|f^{\epsilon}(t,y)|}{\rho _{2}(y)}%
\right) ^{p}\right] \right) <\infty .
\end{equation*}

\noindent Then we can find finite constants $c_{1},c_{2}$ and $c_{3}$ which
do not depend on $\epsilon $, such that for any $q>8$ and any $l,r,v,w>0$
with $18r+9v<q-8,$ there exists $\epsilon _{0}=\epsilon _{0}(t_{0},q,r,v,l)$
such that if $\epsilon \leq \epsilon _{0}<1$ and $\delta w^{-1}>\max
(q/2-r+v/2,(\kappa -n+\alpha )/4\alpha )$ we have%
\begin{eqnarray*}
P\Bigg(\int_{0}^{t_{0}}\left( Y^{\epsilon }(t)\right) ^{2}dt<\epsilon
^{qw},\int_{0}^{t_{0}}\Bigg(\left\vert a^{\epsilon }(t)-\int_{|y|<\epsilon
^{z}}f^{\epsilon }(t,y)G(dy)\right\vert ^{2}+|u^{\epsilon }(t)|^{2}\Bigg)dt
&\geq &l\epsilon ^{w}\Bigg) \\
\leq c_{1}\epsilon ^{rwp}+c_{2}\epsilon ^{wp/4}+c_{3}\exp \left( -\epsilon
^{-vw/2}\right) . &&
\end{eqnarray*}

\noindent Moreover, we have $\epsilon _{0}(t_{0},q,r,v,l)=t_{0}^{-k}\epsilon
_{0}(q,r,v,l)$ for some $k>0.$
\end{lemma}

\begin{proof}
Let $0<C<\infty $ denote a generic constant which varies from line to line
and which does not depend on $\epsilon .$ We begin with some preliminary
remarks. Firstly, the hypotheses of the theorem are sufficient to imply (by
Theorem A6 of \cite{bichteler}) that 
\begin{equation*}
\noindent \sup_{\epsilon >0}\left( \max \left( E\left[ \sup_{0\leq t\leq
t_{0}}|Y^{\epsilon }(t)|^{p}\right] ,E\left[ \sup_{0\leq t\leq
t_{0}}|a^{\epsilon }(t)|^{p}\right] \right) \right) <\infty .
\end{equation*}

\noindent Secondly, by hypothesis we can find previsible processes $%
D_{t}^{\zeta ,\epsilon }$ and $D_{t}^{f,\epsilon }$ and functions $\phi
^{\zeta }$ and $\phi ^{f}$ not depending on $\epsilon $ such that 
\begin{equation}
|\zeta ^{\epsilon }(t,y)|\leq D_{t}^{\zeta ,\epsilon }\phi ^{\zeta }(y)\text{
\ and }|f^{\epsilon }(t,y)|\leq D_{t}^{f,\epsilon }\phi ^{f}(y).\text{\ }
\label{vfc}
\end{equation}

\noindent \noindent Let $D_{t}^{\epsilon }=\max (D_{t}^{\zeta ,\epsilon },$ $%
D_{t}^{f,\epsilon }$ ) and $\phi (y)=\max (\phi ^{\zeta }(y)$ ,$\phi
^{f}(y)) $ and (using the notation of ( \ref{vfconstant}) ) $K=\max
(K_{\zeta },K_{f},1)$, then for some $\epsilon ^{\ast }>0$ we have 
\begin{equation}
\phi (y)\leq K|y|^{\kappa -n+\alpha }  \label{vfc3}
\end{equation}

\noindent for $|y|\leq \epsilon ^{\ast }.$ Consequently taking $\epsilon
\leq \min (\epsilon ^{\ast },1)$ and using the definition of $z$ we see that
for $|y|<\epsilon ^{z}$%
\begin{equation}
\phi (y)\leq K\epsilon ^{z(\kappa -n+\alpha )}=K\epsilon ^{3\delta }.
\label{vfc2}
\end{equation}

\noindent Now, we define 
\begin{equation*}
A=\left\{ \int_{0}^{t_{0}}\left( Y^{\epsilon }(t)\right) ^{2}dt<\epsilon
^{qw},\int_{0}^{t_{0}}\Bigg(\left\vert a^{\epsilon }(t)-\int_{|y|<\epsilon
^{z}}f^{\epsilon }(t,y)G(dy)\right\vert ^{2}+|u^{\epsilon }(t)|^{2}\Bigg)%
dt\geq l\epsilon ^{w}\right\}
\end{equation*}

\noindent and let 
\begin{equation*}
\theta _{t}=|\beta ^{\epsilon }(t)|+|\gamma ^{\epsilon }(t)|+|a^{\epsilon
}(t)|+|u^{\epsilon }(t)|+\int_{|y|<\epsilon ^{z}}(|\zeta ^{\epsilon
}(t,y)|^{2}+|f^{\epsilon }(t,y)|^{2})G(dy).
\end{equation*}

\noindent Taking $\psi =\alpha (\kappa -n+\alpha )^{-1}\leq 1$ we see using (%
\ref{vfc}) and (\ref{vfc2}) that on the set \noindent $\left\{ \sup_{0\leq
t\leq t_{0}}|D_{t}^{\epsilon }|\leq K^{-1}\epsilon ^{-\psi \delta }\right\} $
we have%
\begin{equation}
\underset{0\leq t\leq t_{0}}{\sup }\max (|\zeta ^{\epsilon
}(t,y)|,|f^{\epsilon }(t,y)|)\leq \epsilon ^{-\psi \delta }\epsilon
^{3\delta }\leq \epsilon ^{2\delta }.  \label{bound}
\end{equation}

\noindent Define the stopping time $T=\min (\inf \left\{ s\geq 0:\sup_{0\leq
u\leq s}\theta _{s}>\epsilon ^{-rw}\right\} ,t_{0})$, let $A_{1}=\left\{
T<t_{0}\right\} ,$ $A_{2}=\left\{ \sup_{0\leq t\leq t_{0}}|D_{t}^{\epsilon
}|>K^{-1}\epsilon ^{-\psi \delta }\right\} ,$ $A_{3}=A\cap A_{1}^{c}\cap
A_{2}^{c}$ and observe that 
\begin{equation*}
P(A)\leq P(A_{1})+P(A_{2})+P(A_{3}).
\end{equation*}

\noindent Using (\ref{epsilonunif}) ,the finiteness of (\ref{unif}) and
Chebyshev's inequality gives 
\begin{equation*}
P(A_{1})\leq \epsilon ^{rwp}E\left[ \underset{0\leq t\leq t_{0}}{\sup }%
\theta _{s}^{p}\right] \leq C\epsilon ^{rwp}\text{ and }P(A_{2})\leq
\epsilon ^{\delta \psi p}E\left[ \underset{0\leq t\leq t_{0}}{\sup }D_{t}^{p}%
\right] \leq C\epsilon ^{\delta \psi p},
\end{equation*}

\noindent while on the set $A_{3}$ the processes $a^{\epsilon }$ and $%
Y^{\epsilon }$ satisfy, by virtue of (\ref{bound}), the SDEs%
\begin{eqnarray*}
da^{\epsilon }(t) &=&\beta ^{\epsilon }(t)dt+\underset{i=1}{\overset{d}{\sum 
}}\gamma _{i}^{\epsilon }(t)dW_{t}^{i}+\int_{|y|<\epsilon ^{z}}\zeta
^{\epsilon }(t,y)1_{\left\{ |\zeta ^{\epsilon }(t,y)|<\epsilon ^{2\delta
}\right\} }(\mu -\nu )(dt,dy),\text{ } \\
dY^{\epsilon }(t) &=&a^{\epsilon }(t)dt+\underset{i=1}{\overset{d}{\sum }}%
u_{i}^{\epsilon }(t)dW_{t}^{i}+\int_{|y|<\epsilon ^{z}}f^{\epsilon
}(t,y)1_{\left\{ |f^{\epsilon }(t,y)|<\epsilon ^{2\delta }\right\} }(\mu
-\nu )(dt,dy),\text{ }
\end{eqnarray*}

\noindent with $a^{\epsilon }(0)=\alpha ,$ $Y^{\epsilon }(0)=y.$ We now
define the following processes%
\begin{eqnarray*}
A_{t} &=&\int_{0}^{t}a^{\epsilon }(s)ds\text{, \ }M_{t}=\underset{i=1}{%
\overset{d}{\sum }}\int_{0}^{t}u_{i}^{\epsilon }(s)dW_{s}^{i},\text{ \ }%
Q_{t}=\underset{i=1}{\overset{d}{\sum }}\int_{0}^{t}A(s)\gamma
_{i}^{\epsilon }(s)dW_{s}^{i},\text{ \ } \\
N_{t} &=&\underset{i=1}{\overset{d}{\sum }}\int_{0}^{t}Y^{\epsilon
}(s-)u_{i}^{\epsilon }(s)dW_{s}^{i},\text{ }P_{t}=\int_{0}^{t}\int_{|y|<%
\epsilon ^{z}}f^{\epsilon }(s,y)1_{\left\{ |f^{\epsilon }(s,y)|<\epsilon
^{2\delta }\right\} }(\mu -\nu )(ds,dy),\text{ \ } \\
L_{t} &=&\int_{0}^{t}\int_{|y|<\epsilon ^{z}}Y^{\epsilon }(s-)f^{\epsilon
}(s,y)1_{\left\{ |f^{\epsilon }(s,y)|<\epsilon ^{2\delta }\right\} }(\mu
-\nu )(ds,dy),\text{ } \\
H_{t} &=&\int_{0}^{t}\int_{|y|<\epsilon ^{z}}A(s)\zeta ^{\epsilon
}(s,y)1_{\left\{ |\zeta ^{\epsilon }(s,y)|<\epsilon ^{2\delta }\right\}
}(\mu -\nu )(ds,dy), \\
\text{ \ }J_{t} &=&\int_{0}^{t}\int_{|y|<\epsilon ^{z}}f^{\epsilon
}(s,y)^{2}1_{\left\{ |f^{\epsilon }(s,y)|<\epsilon ^{2\delta }\right\} }(\mu
-\nu )(ds,dy),
\end{eqnarray*}

\noindent and for $\delta _{j}>0,\rho _{j}>0,$ $j\in \left\{ 1,...,7\right\} 
$ define the sets%
\begin{eqnarray*}
B_{1} &=&\left\{ \left\langle N\right\rangle _{T}<\rho _{1},\underset{0\leq
t\leq T}{\sup }|N_{t}|\geq \delta _{1}\right\} ,\text{ }B_{2}=\left\{
\left\langle M\right\rangle _{T}<\rho _{2},\underset{0\leq t\leq T}{\sup }%
|M_{t}|\geq \delta _{2}\right\} , \\
B_{3} &=&\left\{ \left\langle Q\right\rangle _{T}<\rho _{3},\underset{0\leq
t\leq T}{\sup }|Q_{t}|\geq \delta _{3}\right\} ,\text{ }C_{1}=\left\{
\left\langle P\right\rangle _{T}<\rho _{4},\underset{0\leq t\leq T}{\sup }%
|P_{t}|\geq \delta _{4}\right\} , \\
C_{2} &=&\left\{ \left\langle L\right\rangle _{T}<\rho _{5},\sup_{0\leq
t\leq T}|L_{t}|\geq \delta _{5}\right\} ,\text{ \ }C_{3}=\left\{
\left\langle N\right\rangle _{T}<\rho _{6},\underset{0\leq t\leq T}{\sup }%
|N_{t}|\geq \delta _{6}\right\} , \\
C_{4} &=&\left\{ \left\langle J\right\rangle _{T}<\rho _{7},\underset{0\leq
t\leq T}{\sup }|J_{t}|\geq \delta _{7}\right\} .
\end{eqnarray*}

\noindent The exponential martingale inequality for continuous
semimartingales gives $P(B_{j})\leq 2e^{-\delta _{j}^{2}/2\rho _{j}}$ for $%
j=1,2,3.$ Since the jumps in $P$ and $J$ are bounded by $\epsilon ^{2\delta
} $ and $\epsilon ^{4\delta }$ respectively, an application of lemma \ref%
{EMI} gives $\ $%
\begin{equation*}
P(C_{1})\leq 2\exp \left( \frac{-\delta _{4}^{2}}{2(\epsilon ^{2\delta
}\delta _{4}+\rho _{4})}\right) \text{ \ and }P(C_{4})\leq 2\exp \left( 
\frac{-\delta _{7}^{2}}{2(\epsilon ^{4\delta }\delta _{7}+\rho _{7})}\right)
.
\end{equation*}

\noindent \noindent For $C_{2}$ and $C_{3}$ we use the fact that $%
\sup_{0\leq t\leq T}|a^{\epsilon }(t)|\in L^{p}$ and $\sup_{0\leq t\leq
T}|Y^{\epsilon }(t)|\in L^{p}$ uniformly in $\epsilon $ to see 
\begin{eqnarray*}
P(C_{2}) &\leq &P\left( \left\langle L\right\rangle _{T}<\rho
_{5},\sup_{0\leq t\leq T}|L_{t}|\geq \delta _{5},\sup_{0\leq t\leq
T}|Y^{\epsilon }(t)|\leq \epsilon ^{-\delta }\right) \\
&&+P\left( \sup_{0\leq t\leq T}|Y^{\epsilon }(t)|>\epsilon ^{-\delta }\right)
\\
&\leq &2\left( \frac{-\delta _{5}^{2}}{2(\epsilon ^{\delta }\delta _{5}+\rho
_{5})}\right) +C\epsilon ^{\delta p},
\end{eqnarray*}

\noindent where the second term comes from Chebyshev's inequality and the
first follows from lemma \ref{EMI} in concert with the observation that, on
the set $\left\{ \sup_{0\leq t\leq T}|Y^{\epsilon }(t)|\leq \epsilon
^{-\delta }\right\} $, we have 
\begin{equation*}
L_{t}=\int_{0}^{t}\int_{|y|<\epsilon ^{z}}Y^{\epsilon }(s-)f^{\epsilon
}(s,y)1_{\left\{ |f^{\epsilon }(s,y)|<\epsilon ^{2\delta },|Y^{\epsilon
}(s-)|\leq \epsilon ^{-\delta }\right\} }(\mu -\nu )(ds,dy)
\end{equation*}

\noindent for $0\leq t\leq T$. Hence, the jumps in $L$ are bounded by $%
\epsilon ^{\delta }$ on this set (the same argument may also be applied to $%
C_{3}$). \ We now show that $A_{3}\subset \left( \cup _{j=1}^{3}B_{j}\right)
\cup \left( \cup _{j=1}^{4}C_{j}\right) $ whence on choosing appropriate
values for $\delta _{j}$ and $\rho _{j}$ the proof shall be complete. To do
this suppose that $\omega \notin \left( \cup _{j=1}^{3}B_{j}\right) \cup
\left( \cup _{j=1}^{4}C_{j}\right) ,$ $T(\omega )=t_{0},$ $%
\int_{0}^{T}Y_{t}^{\epsilon }(\omega )^{2}dt<\epsilon ^{qw}$ and $%
\sup_{0\leq t\leq T}|D_{t}^{\epsilon }(\omega )|<K^{-1}\epsilon ^{-\psi
\delta }.$ \ Then 
\begin{equation*}
\left\langle N\right\rangle _{T}=\int_{0}^{T}(Y^{\epsilon
}(t-))^{2}|u^{\epsilon }(t)|^{2}dt<\epsilon ^{(-2r+q)w}=:\rho _{1},
\end{equation*}

\noindent and since $\omega \notin B_{1},$ $\sup_{0\leq t\leq T}\left\vert
\sum_{i=1}^{d}\int_{0}^{t}Y^{\epsilon }(s-)u_{i}^{\epsilon
}(s)dW_{s}^{i}\right\vert <\delta _{1}:=\epsilon ^{q_{1}}$, where $%
q_{1}=(q/2-r-v/2)w$. By the same reasoning we have 
\begin{equation*}
\left\langle L\right\rangle _{T}=\int_{0}^{T}\int_{|y|<\epsilon
^{z}}Y^{\epsilon }(t-)^{2}f^{\epsilon }(s,y)^{2}1_{\left\{ |f^{\epsilon
}(t,y)|<\epsilon ^{2\delta }\right\} }G(dy)dt<\epsilon ^{(-2r+q)w}=:\rho
_{5},
\end{equation*}

\noindent since $\omega \notin C_{2}$ we may let $\delta _{5}=\epsilon
^{q_{1}}$ to give $\sup_{0\leq t\leq T}|L_{t}|<\delta _{5}.$ Since we also
have%
\begin{equation*}
\sup_{0\leq t\leq T}\left\vert \int_{0}^{t}Y^{\epsilon }(s-)a^{\epsilon
}(s)ds\right\vert \leq \left( t_{0}\int_{0}^{T}Y^{\epsilon
}(s-)^{2}a^{\epsilon }(s)^{2}ds\right) ^{1/2}<t_{0}^{1/2}\epsilon
^{(-r+q/2)w},
\end{equation*}

\noindent it follows that%
\begin{equation*}
\sup_{0\leq t\leq T}\left\vert \int_{0}^{t}Y^{\epsilon }(s-)dY^{\epsilon
}(s)\right\vert <t_{0}^{1/2}\epsilon ^{(-r+q/2)w}+2\epsilon ^{q_{1}}.
\end{equation*}

\noindent It\^{o}'s formula now gives $Y^{\epsilon
}(t)^{2}=y^{2}+2\int_{0}^{t}Y^{\epsilon }(s-)dY^{\epsilon }(s)+\left\langle
M\right\rangle _{t}+\left[ P\right] _{t},$ and we notice that because 
\begin{eqnarray*}
\left\langle J\right\rangle _{T} &=&\int_{0}^{T}\int_{|y|<\epsilon
^{z}}f^{\epsilon }(s,y)^{4}1_{\left\{ |f^{\epsilon }(s,y)|<\epsilon
^{2\delta }\right\} }G(dy)dt \\
&\leq &\epsilon ^{4\delta }\int_{0}^{T}\int_{|y|<\epsilon ^{z}}f^{\epsilon
}(s,y)^{2}1_{\left\{ |f^{\epsilon }(s,y)|<\epsilon ^{2\delta }\right\}
}G(dy)dt\leq \epsilon ^{4\delta -rw}=:\rho _{7},
\end{eqnarray*}

\noindent and since $\omega \notin C_{4}$ we must have $\sup_{0\leq t\leq
T}|J_{t}|=\sup_{0\leq t\leq T}|\left[ P\right] _{t}-\left\langle
P\right\rangle _{t}|\leq \delta _{7}:=\epsilon ^{2\delta -(r+v)w}.$
Consequently, 
\begin{equation*}
\left\langle M\right\rangle _{t}+\left\langle P\right\rangle _{t}\leq
Y^{\epsilon }(t)^{2}-y^{2}-2\int_{0}^{t}Y^{\epsilon }(s-)dY^{\epsilon
}(s)+\sup_{0\leq t\leq T}|\left[ P\right] _{t}-\left\langle P\right\rangle
_{t}|
\end{equation*}

\noindent and hence,%
\begin{equation*}
\int_{0}^{T}\left\langle M\right\rangle _{t}dt+\int_{0}^{T}\left\langle
P\right\rangle _{t}dt<\epsilon ^{qw}+t_{0}^{3/2}\epsilon
^{(-r+q/2)w}+2t_{0}\epsilon ^{q_{1}}+t_{0}\epsilon ^{2\delta -(r+v)w}.
\end{equation*}

\noindent We notice that $2\delta -(r+v)w>(q-3r)w>q_{1},$ $qw>q_{1}$ and $%
(q/2-r)w>q_{1}$ and so provided 
\begin{equation*}
\epsilon <\min \left( 1,t_{0}^{-1/\left( 2\delta -(r+v)w-q_{1}\right)
},t_{0}^{-3/2\left( \left( -r+q/2\right) w-q_{1}\right) }\right)
\end{equation*}
we get%
\begin{equation*}
\int_{0}^{T}\left\langle M\right\rangle _{t}dt+\int_{0}^{T}\left\langle
P\right\rangle _{t}dt<(2t_{0}+3)\epsilon ^{q_{1}}.
\end{equation*}

\noindent $\left\langle M\right\rangle _{t}$ and $\left\langle
P\right\rangle _{t}$ are increasing processes, so for any $0<\gamma <T$%
\begin{equation*}
\gamma \left\langle M\right\rangle _{T-\gamma }<(2t_{0}+3)\epsilon ^{q_{1}}%
\text{ \ and \ }\gamma \left\langle P\right\rangle _{T-\gamma
}<(2t_{0}+3)\epsilon ^{q_{1}}.
\end{equation*}

\noindent Since these processes are also continuous we get $\left\langle
M\right\rangle _{T}\leq \gamma ^{-1}(2t_{0}+3)\epsilon ^{q_{1}}+\gamma
\epsilon ^{-2rw}$ and $\left\langle P\right\rangle _{T}\leq \gamma
^{-1}(2t_{0}+3)\epsilon ^{q_{1}}+\gamma \epsilon ^{-2rw}$ . By defining $%
\rho _{2}=\rho _{4}:=2(2t_{0}+3)^{1/2}\epsilon ^{-2rw+q_{1}/2}$ and $\gamma
=(2t_{0}+3)^{1/2}\epsilon ^{q_{1}/2}$, we get $\left\langle M\right\rangle
_{T}<\rho _{2}$ and $\left\langle P\right\rangle _{T}<\rho _{4}$, and since $%
\omega \notin B_{2}\cup C_{1}$ we have 
\begin{equation*}
\sup_{0\leq t\leq T}|M_{t}|<\delta _{2}:=\epsilon
^{(q/8-5r/4-5v/8)w}=:\epsilon ^{q_{2}},\text{ \ }\sup_{0\leq t\leq
T}|P_{t}|<\delta _{4}=\epsilon ^{q_{2}}.
\end{equation*}

\noindent Since $\int_{0}^{T}Y^{\epsilon }(t)^{2}dt<\epsilon ^{qw}$
Chebyshev's inequality gives%
\begin{equation*}
Leb\left\{ t\in \lbrack 0,T]:\left\vert Y_{t}^{\epsilon }(\omega
)\right\vert \geq \epsilon ^{qw/3}\right\} \leq \epsilon ^{qw/3}
\end{equation*}

\noindent so that 
\begin{equation*}
Leb\left\{ t\in \lbrack 0,T]:\left\vert y+A_{t}(\omega )\right\vert \geq
\epsilon ^{qw/3}+2\epsilon ^{q_{2}}\right\} \leq \epsilon ^{qw/3}.
\end{equation*}

\noindent Then, for each $t\in \lbrack 0,T],$ there exists some $s\in
\lbrack 0,T]$ such that $|s-t|\leq \epsilon ^{qw/3}$ and $|y+A_{s}(\omega
)|<\epsilon ^{qw/3}+2\epsilon ^{q_{2}}$, which yields%
\begin{equation*}
|y+A_{t}|\leq |y+A_{s}|+\left\vert \int_{s}^{t}a^{\epsilon }(\tau )d\tau
\right\vert <(1+\epsilon ^{-rw})\epsilon ^{qw/3}+2\epsilon ^{q_{2}}.
\end{equation*}

\noindent In particular we have $|y|<(1+\epsilon ^{-rw})\epsilon
^{qw/3}+2\epsilon ^{q_{2}}$ and, for all $t\in \lbrack 0,T]$, since $%
q_{2}<(q/3-r)w$, we have%
\begin{equation*}
|A_{t}|<2\left( (1+\epsilon ^{-rw})\epsilon ^{qw/3}+2\epsilon
^{q_{2}}\right) \leq 8\epsilon ^{q_{2}}.
\end{equation*}

\noindent This implies that%
\begin{eqnarray*}
\left\langle Q\right\rangle _{T} &=&\int_{0}^{T}A(t)^{2}|\gamma ^{\epsilon
}(t)|^{2}dt<64t_{0}\epsilon ^{2q_{2}-2rw}=:\rho _{3} \\
\left\langle H\right\rangle _{T} &=&\int_{0}^{T}\int_{|y|<\epsilon
^{z}}A(t)^{2}\zeta ^{\epsilon }(t,y)^{2}1_{\left\{ |\zeta ^{\epsilon
}(t,y)|<\epsilon ^{2\delta }\right\} }G(dy)dt\leq \rho _{3}=:\rho _{6},
\end{eqnarray*}

\noindent and since $\omega \notin B_{3}\cup C_{3}$ we must have 
\begin{eqnarray*}
\underset{0\leq t\leq T}{\sup }|Q_{t}| &=&\underset{0\leq t\leq T}{\sup }%
\left\vert \overset{d}{\underset{i=1}{\sum }}\int_{0}^{t}A(s)\gamma
_{i}^{\epsilon }(s)dW_{i}(s)\right\vert <\delta _{3}:=\epsilon
^{(q/8-9r/4-9v/8)w}=:\epsilon ^{q_{3}} \\
\underset{0\leq t\leq T}{\sup }|H_{t}| &=&\underset{0\leq t\leq T}{\sup }%
\left\vert \int_{0}^{t}\int_{|y|<\epsilon ^{z}}A(s)\zeta ^{\epsilon
}(s,y)1_{\left\{ |\zeta ^{\epsilon }(s,y)|<\epsilon ^{2\delta }\right\}
}(\mu -\nu )(ds,dy)\right\vert <\delta _{6}:=\epsilon ^{q_{3}}.
\end{eqnarray*}

\noindent Now we observe using (\ref{vfc}),(\ref{vfc3}), condition \ref%
{condition2} , $\sup_{0\leq t\leq T}|D_{t}^{\epsilon }(\omega
)|<K^{-1}\epsilon ^{-\psi \delta }$, the definition of $\psi $, and the fact
that $\phi ^{f}$ does not depend on $\epsilon $ 
\begin{eqnarray*}
\int_{0}^{t_{0}}\left\vert \int_{|y|<\epsilon ^{z}}f^{\epsilon
}(t,y)G(dy)\right\vert ^{2}dt &\leq &t_{0}\left( \epsilon ^{-\delta \psi
}\int_{|y|<\epsilon ^{z}}|y|^{\kappa -n+\alpha }G(dy)\right) ^{2} \\
&\leq &Ct_{0}\epsilon ^{-2\delta \psi +2z\alpha }=Ct_{0}\epsilon ^{4\delta
\alpha /(\kappa -n+\alpha )}.
\end{eqnarray*}

\noindent An application of \noindent It\^{o}'s formula then gives 
\begin{eqnarray*}
&&\int_{0}^{T}\left( \left\vert a^{\epsilon }(t)-\int_{|y|<\epsilon
^{z}}f^{\epsilon }(t,y)G(dy)\right\vert ^{2}+|u^{\epsilon }(t)|^{2}\right) dt
\\
&\leq &2\int_{0}^{T}a^{\epsilon }(t)^{2}dt+\int_{0}^{T}|u^{\epsilon
}(t)|^{2}dt+2\int_{0}^{T}\left\vert \int_{|y|<\epsilon ^{z}}f^{\epsilon
}(t,y)G(dy)\right\vert ^{2}dt \\
&\leq &2\int_{0}^{T}a^{\epsilon }(t)dA(t)+\left\langle M\right\rangle
_{T}+2Ct_{0}\epsilon ^{4\delta \alpha /(\kappa -n+\alpha )} \\
&=&2\Bigg(a^{\epsilon }(T)A(T)-\int_{0}^{T}A(t)\beta ^{\epsilon }(t)dt-%
\underset{i=1}{\overset{d}{\sum }}\int_{0}^{T}A(t)\gamma _{i}^{\epsilon
}(t)dW_{t}^{i} \\
&&-\int_{0}^{T}\int_{|y|<\epsilon ^{z}}A(s)\zeta ^{\epsilon }(s,y)1_{\left\{
|\zeta ^{\epsilon }(s,y)|<\epsilon ^{2\delta }\right\} }(\mu -\nu )(ds,dy)%
\Bigg)+\left\langle M\right\rangle _{T} \\
&&+2Ct_{0}\epsilon ^{4\delta \alpha /(\kappa -n+\alpha )} \\
&\leq &16(1+t_{0})\epsilon ^{q_{2}-rw}+4\epsilon
^{q_{3}}+4(2t_{0}+3)^{1/2}\epsilon ^{-2rw+q_{1}/2}+2Ct_{0}\epsilon ^{4\delta
\alpha /(\kappa -n+\alpha )} \\
&\leq &l\epsilon ^{w}
\end{eqnarray*}

\noindent provided%
\begin{multline}
\epsilon <\min \Bigg(\left( \frac{l}{16(1+t_{0})}\right)
^{(q_{2}-rw-w)^{-1}},\left( \frac{l}{4}\right) ^{(q_{3}-w)^{-1}},  \label{a}
\\
\left( \frac{l}{4\left( 2t_{0}+3\right) ^{1/2}}\right)
^{(-2rw+q_{1}/2-w)^{-1}},\left( \frac{l}{2Ct_{0}}\right) ^{\left( \frac{%
4\delta \alpha }{\kappa -n+\alpha }-w\right) ^{-1}}\Bigg).  \notag
\end{multline}

\noindent Where the last inequality follows from $q_{2}-rw>w,$ $q_{3}>w,$ $%
q_{1}/2-2rw>w$ and $\delta >w(\kappa -n+\alpha )/4\alpha $. \ Finally, by
the choice of $\delta _{j}$ and $\rho _{j}$ and the assumption that $\delta
>(-r+q/2+v/2)w$ (which also implies that $\delta >-rw+q_{1}/2-q_{2}/4$ and $%
\delta >2q_{2}-2rw-q_{3}$) we see that $\epsilon ^{2\delta }\delta _{4}<\rho
_{4},$ $\epsilon ^{\delta }\delta _{5}<\rho _{5},$ $\epsilon ^{\delta
}\delta _{6}<\rho _{6}$ and $\epsilon ^{4\delta }\delta _{7}<\rho _{7}.$
Therefore this choice for $\delta _{j}$ and $\rho _{j}$ enable us to deduce
that 
\begin{eqnarray*}
P\left( \underset{j=1}{\overset{3}{\cup }}B_{j}\right) &\leq &2\Bigg(\exp
\left( -\frac{1}{2}\epsilon ^{-vw}\right) +\exp \left( -\frac{1}{%
4(2t_{0}+3)^{1/2}}\epsilon ^{-vw}\right) \\
&&+\exp \left( -\frac{1}{128t_{0}}\epsilon ^{-vw}\right) \Bigg),
\end{eqnarray*}

\noindent and%
\begin{eqnarray*}
P\left( \underset{j=1}{\overset{4}{\cup }}C_{j}\right) &\leq &2\Bigg(2\exp
\left( -\frac{1}{4}\epsilon ^{-vw}\right) +\exp \left( -\frac{1}{%
8(2t_{0}+3)^{1/2}}\epsilon ^{-vw}\right) \\
&&+\exp \left( -\frac{1}{256t_{0}}\epsilon ^{-vw}\right) +C\epsilon ^{\psi
\delta p}\Bigg).
\end{eqnarray*}

\noindent The proof is finished on noting that $\delta \psi >w/4,$ and the
dependence of $\epsilon _{0}$ on $t_{0}$ follows immediately from the proof.
\end{proof}

\section{\protect\bigskip Uniform H\"{o}rmander condition}

\bigskip We now present our uniform H\"{o}rmander condition.

\begin{condition}[UH]
\label{UH}Let $V_{0}=Z-\frac{1}{2}\sum_{i=1}^{d}D V_{i}V_{i}$ and assume
that condition \ref{condition1} holds. Recursively define the following
families of vector fields%
\begin{gather*}
\mathcal{L}_{0}=\left\{ V_{1},...,V_{d}\right\} \\
\mathcal{L}_{n+1}=\mathcal{L}_{n}\cup \left\{ \lbrack V_{i},K],\text{ }%
i=1,...,d:\text{ }K\in \mathcal{L}_{n}\right\} \\
\cup \left\{ \left[ V_{0},K\right] -\int_{E}[Y,K](\cdot ,y)G(dy):K\in 
\mathcal{L}_{n}\right\} .
\end{gather*}

\noindent Then there exists some smallest integer $j_{0}\geq 1$ and a
constant $c>0$ such that for any $u\in 
%TCIMACRO{\U{211d} }%
%BeginExpansion
\mathbb{R}
%EndExpansion
^{e}$ with $|u|=1$ we have%
\begin{equation*}
\underset{x\in 
%TCIMACRO{\U{211d} }%
%BeginExpansion
\mathbb{R}
%EndExpansion
^{e}}{\inf }\overset{j_{0}}{\underset{j=0}{\sum }}\underset{K\in \mathcal{L}%
_{j}}{\sum }\left( u^{T}K(x)\right) ^{2}\geq c
\end{equation*}
\end{condition}

\noindent The next important result is a development of an idea presented in 
\cite{teichmann}, it enables us to estimate the Malliavin covariance matrix
on a time interval where the Poisson random measure records no jumps of size
greater than some truncation parameter. As in \cite{teichmann} the key idea
is to make explicit the dependence of the estimate on the length of the time
interval under consideration.

\begin{theorem}
\label{cthm}Let $t>0$ and let $x_{t}$ satisfy the SDE 
\begin{equation*}
x_{t}=x+\int_{0}^{t}Z(x_{s-})ds+\int_{0}^{t}V(x_{s-})dW_{s}+\int_{0}^{t}%
\int_{E}Y(y,x_{s-})(\mu -\nu )(dy,ds)
\end{equation*}%
\noindent and assume that the following conditions are satisfied :%
\begin{equation}
Z,V_{1},..,V_{d}\in C_{b}^{\infty }(%
%TCIMACRO{\U{211d} }%
%BeginExpansion
\mathbb{R}
%EndExpansion
^{e}),
\end{equation}

\noindent for every $y\in E$ $Y(\cdot ,y)\in C_{b}^{\infty }(%
%TCIMACRO{\U{211d} }%
%BeginExpansion
\mathbb{R}
%EndExpansion
^{e})$ and, for some $\rho _{2}\in L^{2,\infty }(G)$ and every $n\in 
%TCIMACRO{\U{2115} }%
%BeginExpansion
\mathbb{N}
%EndExpansion
\cup \left\{ 0\right\} $%
\begin{equation}
\underset{y\in E,x\in 
%TCIMACRO{\U{211d} }%
%BeginExpansion
\mathbb{R}
%EndExpansion
^{e}}{\sup }\frac{1}{\rho _{2}(y)}|D_{1}^{n}Y(x,y)|<\infty ,
\label{Lpinequalities}
\end{equation}%
\begin{equation*}
\underset{y\in E,x\in 
%TCIMACRO{\U{211d} }%
%BeginExpansion
\mathbb{R}
%EndExpansion
^{e}}{\sup }|\left( I+D_{1}Y(x,y)\right) ^{-1}|<\infty \text{ \ and }%
\underset{x\in 
%TCIMACRO{\U{211d} }%
%BeginExpansion
\mathbb{R}
%EndExpansion
^{e}}{\sup }|\left( I+D_{1}Y(x,\cdot )\right) ^{-1}|\in L^{\infty }(G).\text{
}
\end{equation*}

\noindent Further assume conditions \ref{condition1}, \ref{condition2} , \ref%
{condition3} \ and condition (UH) hold. For some $0<t<t_{0},$ $\delta
,\alpha >0$ and $z=3\delta (\kappa -n+\alpha )^{-1}$ define the set $%
A_{t}=A_{t}(\epsilon )$ by 
\begin{equation*}
A_{t}=\left\{ \omega :\text{ }\left( \text{supp }\mu (\cdot ,\cdot )\right)
\cap \lbrack 0,t)\times E\subseteq \lbrack 0,t)\times \left\{ |y|\leq
\epsilon ^{z}\right\} \right\} .
\end{equation*}

\noindent Then, $P\left( \left\{ \sup_{0\leq s\leq t}|x_{s}-x_{s}(\epsilon
)|>0\right\} \cap A_{t}\right) =0$, where $x_{t}(\epsilon )$ is the solution
to the SDE%
\begin{eqnarray}
dx_{t}(\epsilon ) &=&\left( Z(x_{t-}(\epsilon ))-\int_{|y|\geq \epsilon
^{z}}Y(x_{t-}(\epsilon ),y)G(dy)\right) dt+V(x_{t-}(\epsilon ))dW_{t}  \notag
\\
&&\text{ \ \ \ \ \ \ \ \ \ \ \ \ \ \ \ \ \ \ \ \ \ \ \ \ \ \ \ \ \ \ \ \ \ \
\ \ \ \ \ }+\int_{|y|<\epsilon ^{z}}Y(x_{t-}(\epsilon ),y)(\mu -\nu )(dy,dt),
\label{conditionedprocess}
\end{eqnarray}
Moreover if we let the reduced Malliavin covariance matrix associated with $%
x_{t}(\epsilon )$ be denoted by $C_{t}(\epsilon )$ then we have for any $%
p\geq 1$ and some $\epsilon _{0}(p)>0,$ $K(p)\geq 1$, that 
\begin{equation*}
\underset{|u|=1}{\sup }P(\left\{ u^{T}C_{t}u\leq \epsilon \right\} \cap
A_{t})=\underset{|u|=1}{\sup }P(u^{T}C_{t}(\epsilon )u\leq \epsilon )\leq
\epsilon ^{p}
\end{equation*}
\noindent for $0\leq \epsilon \leq t^{K(p)}\epsilon _{0}(p),$ provided that 
\begin{equation*}
16 \delta > \max \left(8-r+\frac{v}{2},\frac{\kappa -n+\alpha }{4\alpha}
\right),
\end{equation*}
where $r,v>0$ are such that $18r + 9v<8$.
\end{theorem}

\begin{proof}
The indistinguishability of the processes $x$ and $x(\epsilon )$ on $A_{t}$
is a trivial. For the remainder of the proof we first note that condition
(UH) uniform enables us to identify a smallest integer $j_{0}$ and a
constant $c>0$ such that, for any $u\in 
%TCIMACRO{\U{211d} }%
%BeginExpansion
\mathbb{R}
%EndExpansion
^{e}$ with $|u|=1$ 
\begin{equation*}
\underset{x\in 
%TCIMACRO{\U{211d} }%
%BeginExpansion
\mathbb{R}
%EndExpansion
^{e}}{\inf }\overset{j_{0}}{\underset{j=0}{\sum }}\underset{K\in \mathcal{L}%
_{j}}{\sum }\left( u^{T}K(x)\right) ^{2}\geq c.
\end{equation*}

\noindent For $j=0,1,...,j_{0}$ set $m(j)=2^{-4j}$ and define 
\begin{equation*}
E_{j}=\left\{ \underset{K\in \mathcal{L}_{j}}{\sum }\int_{0}^{t}\left(
u^{T}(\epsilon )J_{0\leftarrow s}(\epsilon )K(x_{s}(\epsilon ))\right)
^{2}ds\leq \epsilon ^{m(j)}\right\} ,
\end{equation*}

\noindent where $J_{t\leftarrow 0}(\epsilon )$ denotes the Jacobian of the
flow associated with $x_{t}(\epsilon )$ and $J_{0\leftarrow t}(\epsilon )$
denotes its inverse (which exists by the assumptions on the vector fields as
in theorem \ref{jacob}). It is straight forward to note, using (\ref%
{Lpinequalities}), $L^{p}$ inequalities for stochastic integrals based on
Poisson random measures (see \cite{bichteler}, lemma A.14) and Gronwall's
inequality that for any $p<\infty $ 
\begin{equation}
\underset{\epsilon \geq 0}{\sup }\mathbb{E}\left[ \underset{0\leq s\leq t}{%
\sup }|J_{t\leftarrow 0}(\epsilon )|^{p}\right] <\infty .  \label{epsLp}
\end{equation}%
Let $C$ denote a constant which varies from line to line and does not depend
on $\epsilon $. Then, as usual we have 
\begin{equation*}
\left\{ u^{T}C_{t}(\epsilon )u\leq \epsilon \right\} =E_{0}\subset \left(
E_{0}\cap E_{1}^{c}\right) \cup (E_{1}\cap E_{2}^{c})\cup ...\cup
(E_{j_{0}-1}\cap E_{j_{0}}^{c})\cup F
\end{equation*}

\noindent where $F=E_{0}\cap E_{1}\cap ...\cap E_{j_{0}}.$ Define the
stopping time%
\begin{equation*}
S=\min \left( \inf \left\{ s\geq 0:\underset{0\leq z\leq s}{\sup }%
|J_{0\leftarrow z}(\epsilon )-I|\geq \frac{1}{2}\right\} ,t\right) ,
\end{equation*}

\noindent and notice that by choosing $0<\beta <m(j_{0})$ we discover that $%
P(F)\leq P(S<\epsilon ^{\beta })\leq C\epsilon ^{q\beta /2}$ for $\epsilon
\leq \epsilon _{1}$and any $q\geq 2$ (see \cite{Nu06} and \cite{teichmann}
for details), where as in \cite{teichmann}, $\epsilon _{1}$ satisfies%
\begin{equation*}
\epsilon _{1}<\min \left( t^{1/\beta },\left( \frac{c}{4(j_{0}+1)}\right)
^{1/(m(j_{0})-\beta )}\right) .
\end{equation*}

\noindent We notice that for any $K\in C_{b}^{\infty }(%
%TCIMACRO{\U{211d} }%
%BeginExpansion
\mathbb{R}
%EndExpansion
^{e})$ we have 
\begin{gather*}
d(u^{T}J_{0\leftarrow t}(\epsilon )K(x_{t}(\epsilon ))=u^{T}J_{0\leftarrow
t-}(\epsilon )\Bigg(\left[ V_{0},K\right] (x_{t-}(\epsilon ))-\int_{E}\left[
Y,K\right] (x_{t-}(\epsilon ),y)G(dy) \\
+\frac{1}{2}\overset{d}{\underset{i=1}{\sum }}\left[ V_{i},\left[ V_{i},K%
\right] \right] (x_{t-}(\epsilon )) \\
+\int_{|y|<\epsilon ^{z}}((I+D_{1}Y(x_{t-}(\epsilon
),y)^{-1})K(x_{t-}(\epsilon )+Y(x_{t-}(\epsilon ),y))-K(x_{t-}(\epsilon
))G(dy))\Bigg)dt \\
+u^{T}J_{0\leftarrow t-}(\epsilon )\overset{d}{\underset{i=1}{\sum }}\left[
V_{i},K\right] (x_{t-}(\epsilon ))dW_{t}^{i} \\
+u^{T}J_{0\leftarrow t-}(\epsilon )\int_{|y|<\epsilon^{z}}
((I+D_{1}Y(x_{t-}(\epsilon ),y)^{-1})K(x_{t-}(\epsilon )+Y(x_{t-}(\epsilon
),y))-K(x_{t-}(\epsilon ))(\mu -\nu )(dy,dt).
\end{gather*}

\noindent We now verify the conditions of lemma \ref{Norris} in the case
where 
\begin{eqnarray*}
Y^{\epsilon }(t) &=&u^{T}J_{0\leftarrow t}(\epsilon )K(x_{t}(\epsilon )) \\
a^{\epsilon }(t) &=&u^{T}J_{0\leftarrow t}(\epsilon )\Bigg(\left[ V_{0},K%
\right] (x_{t}(\epsilon ))-\int_{E}\left[ Y,K\right] (x_{t}(\epsilon
),y)G(dy) \\
&&+\frac{1}{2}\overset{d}{\underset{i=1}{\sum }}\left[ V_{i},\left[ V_{i},K%
\right] \right] (x_{t}(\epsilon )) \\
&&+\int_{|y|<\epsilon ^{z}}((I+D_{1}Y(x_{t}(\epsilon
),y)^{-1})K(x_{t}(\epsilon )+Y(x_{t}(\epsilon ),y))-K(x_{t}(\epsilon ))G(dy))%
\Bigg). \\
&=:&u^{T}J_{0\leftarrow t}\tilde{K}(x_{t}(\epsilon)),
\end{eqnarray*}
where $\tilde{K}\in C_{b}^{\infty}(\mathbb{R}^{e})$. \noindent To do this we
observe, using the notation of lemma \ref{Norris} that%
\begin{equation*}
f^{\epsilon}(t,y)=u^{T}J_{0\leftarrow t-}(\epsilon
)((I+D_{1}Y(x_{t-}(\epsilon ),y)^{-1})K(x_{t-}(\epsilon )+Y(x_{t-}(\epsilon
),y))-K(x_{t-}(\epsilon ))
\end{equation*}

\noindent and hence for some $0<C<\infty $ 
\begin{eqnarray*}
|f^{\epsilon}(t,y)| &\leq &C\left\vert u^{T}J_{0\leftarrow t-}(\epsilon )
\right\vert \max\left(\underset{x\in \mathbb{R}^{e}}{\sup}\left\vert K(x
)\right\vert ,\underset{x\in \mathbb{R}^{e}}{\sup}\left\vert DK(x
)\right\vert\right) \\
&& \bigg(\underset{x \in \mathbb{R}^{e}, y\in E}{\sup }\left\vert \left(
I+D_{1}Y(x ,y)\right)^{-1}\right\vert |D_{1}Y(x_{t-}(\epsilon ),y)|
+|Y(x_{t-}(\epsilon ),y)|\bigg).
\end{eqnarray*}

\noindent Condition \ref{condition3} then gives that $|f^{\epsilon}(t,y)|%
\leq C\left\vert u^{T}J_{0\leftarrow t-}(\epsilon )\right\vert \phi (y)$
where $\phi \in L_{+}^{1}(G)$ does not depend on $\epsilon $, $C=C(K)<\infty 
$ and where and for some $\alpha >0$ (which does not depend on $\epsilon $
or $K$!) we have%
\begin{equation*}
\underset{y\rightarrow 0}{\lim \sup }\frac{\phi (y)}{|y|^{\kappa -n+\alpha }}%
\text{ }<\infty .
\end{equation*}
Finally, using the notation of (\ref{epsilonunif}), we notice that
Cauchy-Schwarz gives 
\begin{equation}
|u^{T}J_{0\leftarrow t-}(\epsilon )|\leq \underset{i=1}{\overset{e}{\sum}}%
|e_{i}^{T}J_{0\leftarrow t-}(\epsilon)|^{2}=:D_{t}^{f,\epsilon},
\label{uniformity}
\end{equation}
where $e_{i}$ is the standard basis in $\mathbb{R}^{e}$. Hence by (\ref%
{epsLp}) we have for any $p<\infty$ 
\begin{align*}
\underset{\epsilon\geq 0}{\sup}E\left[\underset{0\leq s\leq t}{\sup}%
\left(D_{s}^{f,\epsilon}\right)^{p}\right]<\infty.
\end{align*}
\noindent We have therefore verified the conditions of lemma \ref{Norris}
for the process $f^{\epsilon}(t,y)$. They may be also checked for the
process $\zeta^{\epsilon}(t,y)$ in the same manner. The other hypotheses of
lemma \ref{Norris} are trivial to verify so we apply this lemma with $%
z=3\delta (\kappa -n+\alpha )^{-1},$ and with $q=16$, $r,v>0$ such that $%
18r+9v<8$ and $16\delta >\max \left( 8-r+v/2,(\kappa -n+\alpha )/4\alpha
\right) $ to deduce that for $j\in \left\{ 0,1,...,j_{0}-1\right\} $%
\begin{align*}
P(E_{j}\cap E_{j_{+1}}^{c})=P\Bigg(\underset{K\in \mathcal{L}_{j}}{\sum }&
\int_{0}^{t}\left( u^{T}J_{0\leftarrow s}(\epsilon )K(x_{s}(\epsilon
))\right) ^{2}ds\leq \epsilon ^{m(j)}, \\
\text{ \ \ \ \ \ \ \ \ \ \ \ \ \ \ \ \ \ \ \ \ \ \ \ \ \ \ \ \ \ \ \ \ }& 
\text{\ }\underset{K\in \mathcal{L}_{j+1}}{\sum }\int_{0}^{t}\left(
u^{T}J_{0\leftarrow s}(\epsilon )K(x_{s}(\epsilon ))\right) ^{2}ds>\epsilon
^{m(j+1)}\Bigg) \\
\leq \underset{K\in \mathcal{L}_{j}}{\sum }P\Bigg(\int_{0}^{t}(v^{T}&
J_{0\leftarrow s}(\epsilon )K(x_{s}(\epsilon ))^{2}ds\leq \epsilon ^{m(j)},
\\
\text{ \ \ \ \ \ \ \ \ \ \ \ \ \ \ \ \ }\text{\ }\underset{k=1}{\overset{d}{%
\sum }}\int_{0}^{t}(u^{T}J_{0\leftarrow s}(\epsilon )& V_{k}(x_{s}(\epsilon
))^{2}ds+\int_{0}^{t}u^{T}J_{0\leftarrow s-}(\epsilon )\Bigg(\left[ V_{0},K%
\right] (x_{s-}(\epsilon )) \\
\text{ \ \ \ \ \ \ \ \ \ \ \ \ \ }-\int_{E}\left[ Y,K\right]
(x_{s-}(\epsilon ),y)& G(dy)+\frac{1}{2}\overset{d}{\underset{i=1}{\sum }}%
\left[ V_{i},\left[ V_{i},V_{k}\right] \right] (x_{s-}(\epsilon ))\Bigg)%
ds\left. >\right. \frac{\epsilon ^{m(j+1)}}{n(j)}\Bigg).
\end{align*}

\noindent Which \noindent is $o(\epsilon ^{p})$ for $\epsilon \leq \epsilon
_{2}(p)$ where $\epsilon _{2}$ can be chosen as $\epsilon _{3}t^{-k^{\ast }}$
for some $k^{\ast }>0$ and where $\epsilon _{3}$ is independent of $t$.
Setting $\epsilon _{0}=\min (\epsilon _{1},\epsilon _{2})$ and noticing by (%
\ref{uniformity}) that all the estimates are uniform over $|u|=1$ gives the
result.
\end{proof}

\section{C$^{\infty }$ density under the H\"{o}rmander condition}

We now state and prove our main result

\begin{theorem}
\bigskip Suppose that $x_{t}$ is the solution to the SDE 
\begin{equation*}
x_{t}=x+\int_{0}^{t}Z(x_{s-})ds+\int_{0}^{t}V(x_{s-})dW_{s}+\int_{0}^{t}%
\int_{E}Y(y,x_{s-})(\mu -\nu )(dy,ds)
\end{equation*}

\noindent and that the conditions of theorem \ref{cthm} are in force. Then,
for any $t_{0}>0$ the law of $x_{t_{0}}$ has a $C^{\infty }$ density with
respect to Lebesgue measure under the uniform H\"{o}rmander condition \ref%
{UH} provided, in the notation of theorem \ref{cthm}, we have 
\begin{equation}
16m(j_{0})>3(\kappa -n)\max \left( \frac{8-r+v/2}{\kappa -n+\alpha },\frac{1%
}{4\alpha }\right) .  \label{bracketcond}
\end{equation}%
\noindent
\end{theorem}

\begin{remark}
Note that (\ref{bracketcond}) is always true when $\kappa =n.$
\end{remark}

\begin{proof}
By Theorem \ref{criterion} it suffices to check that $\left\vert
C_{t_{0}}^{-1}\right\vert \in L^{p}$ for all $p\geq 2.$ Let $\Lambda= 
\underset{|u|=1}{\inf}u^{T}C_{t_{0}}u$ be the smallest eigenvalue of $%
C_{t_{0}}$. Then it is sufficient to show that $\Lambda^{-1} \in L^{p}$ for
all $p\geq 2$. However, we may write 
\begin{eqnarray}
E[\Lambda ^{-p}]=C_{1} \int_{0}^{\infty }\epsilon ^{-k}P(\Lambda \leq
\epsilon^{2} )d\epsilon \leq C_{2}+C_{3}\int_{0}^{1}\epsilon ^{-k}P(\Lambda
\leq\epsilon^{2} )d\epsilon ,  \notag  \label{inverse}
\end{eqnarray}
for some $k>1$. By a routine compactness argument we may show (see \cite%
{norris}) that 
\begin{align*}
P(\Lambda \leq \epsilon) \leq C_{2}\epsilon^{-e}\underset{|u|=1}{\sup}%
P(u^{T}C_{t_{0}}u\leq \epsilon),
\end{align*}
so that for some $k^{\prime }>1$ 
\begin{equation}
E[\Lambda ^{-p}]\leq C_{3}+C_{4}\int_{0}^{1} \epsilon ^{-k^{\prime }}%
\underset{|u|=1}{\sup}P(u^{T}C_{t_{0}}u\leq \epsilon^{2})d\epsilon .
\label{lambda}
\end{equation}
Now we define a Poisson process $N_{\epsilon }$ on $%
%TCIMACRO{\U{211d} }%
%BeginExpansion
\mathbb{R}
%EndExpansion
^{+}$ for $\epsilon >0$ by 
\begin{equation*}
N_{\epsilon }(t)=\int_{0}^{t}\int_{|y|>\epsilon ^{z}}\mu (dy,ds),
\end{equation*}

\noindent whose rate is given as%
\begin{equation*}
\lambda (\epsilon )=\int_{|y|>\epsilon ^{z}}G(dy).
\end{equation*}%
By (\ref{cond2a}) we know that%
\begin{equation}
\underset{\epsilon \rightarrow 0}{\lim \sup }\frac{\lambda (\epsilon )}{%
f\left( \epsilon \right) }<\infty  \label{rategrowth}
\end{equation}%
We may find a (random) subinterval $[t_{1},t_{2})\subseteq \lbrack 0,t_{0})$
such that $t_{2}-t_{1}\geq t_{0}(N_{\epsilon }(t_{0})+1)^{-1}$ on which the
Poisson random measure $\mu $ records no jumps of absolute value greater
than $\epsilon ^{z}$ and, as such, the underlying process $x_{t}$ solves the
SDE (\ref{conditionedprocess}) started at $x_{t_{1}}$ on this interval. We
emphasize the dependence of $C_{t_{0}}$on the starting point $(x,I)$ of the
process $(x_{t},J_{0\leftarrow t})$. Then, using the fact that $%
J_{0\leftarrow t}^{x,V}=VJ_{0\leftarrow t}^{x,I},$ $J_{0\leftarrow
t}=J_{t\leftarrow 0}^{-1}\ $ , the (strong) Markov property, and the two
observations that $t_{2}-t_{1}\geq t_{0}(N_{\epsilon }(t_{0})+1)^{-1}$ and 
\begin{equation*}
span\{u^{T}J_{0\leftarrow t}^{x,I}:u\in \mathbb{R}^{e},|u|=1\}=\mathbb{R}%
^{e}\,\,\,\,\text{a.s. for every $t>0$ and $x\in \mathbb{R}^{e}$}
\end{equation*}%
we see that for any $q<\infty $ 
\begin{eqnarray}
\underset{|u|=1}{\sup }P(u^{T}C_{t_{0}}^{x,I}u\leq \epsilon ^{2}) &\leq &%
\underset{|u|=1}{\sup }P\left(
u^{T}C_{t_{1},t_{2}}^{x_{t_{1}}^{x},J_{0\leftarrow t_{1}}^{x,I}}u\leq
\epsilon ^{2}\right)  \notag \\
&=&\underset{|u|=1}{\sup }P\left( u^{T}J_{0\leftarrow
t_{1}}^{x,I}C_{t_{1},t_{2}}^{x_{t_{1}}^{x},I}\left( J_{0\leftarrow
t_{1}}^{x,I}\right) ^{T}u\leq \epsilon ^{2}\right)  \notag \\
&=&\underset{|u|=1}{\sup }P\left( \frac{u^{T}J_{0\leftarrow
t_{1}}^{x,I}C_{t_{1},t_{2}}^{x_{t_{1}}^{x},I}\left( J_{0\leftarrow
t_{1}}^{x,I}\right) ^{T}u}{|u^{T}J_{0\leftarrow t_{1}}^{x,I}|^{2}}\leq \frac{%
\epsilon ^{2}}{{|u^{T}J_{0\leftarrow t_{1}}^{x,I}|^{2}}}\right)  \notag \\
&\leq &\underset{|u|=1}{\sup }P\left(
u^{T}C_{t_{1},t_{2}}^{x_{t_{1}}^{x},I}u\leq \epsilon \right) +\underset{|u|=1%
}{\sup }P\left( |u^{T}J_{0\leftarrow t_{1}}^{x,I}|^{-1}\geq \epsilon
^{-1/2}\right)  \notag \\
&=&\underset{|u|=1}{\sup }P\left(
u^{T}C_{t_{2}-t_{1}}^{x_{t_{1}}^{x},I}(\epsilon )u\leq \epsilon \right)
+O(\epsilon ^{q})  \notag \\
&\leq &\underset{|u|=1}{\sup }P\left( u^{T}C_{t_{0}(N_{\epsilon
}(t_{0})+1)^{-1}}^{x_{t_{1}},I}(\epsilon )u\leq \epsilon \right) +O(\epsilon
^{q}).  \label{chain}
\end{eqnarray}%
An application of theorem \ref{cthm} yields%
\begin{equation*}
\underset{|u|=1}{\sup }P\left( u^{T}C_{t_{0}(N_{\epsilon
}(t_{0})+1)^{-1}}^{x_{t_{1}},I}(\epsilon )u\leq \epsilon \right) \,\,\,\text{%
is}\,\,\,O(\epsilon ^{q})
\end{equation*}%
for any $q\geq 2$ \noindent if $\epsilon \leq \epsilon
_{0}t_{0}^{1/K(q)}(N_{\epsilon }(t_{0})+1)^{-1/K(q)}$ provided that \newline
$\delta >\max \left( 8-r+v/2,(\kappa -n+\alpha )/4\alpha \right) .$ From
this, (\ref{lambda}) and (\ref{chain}) we get that 
\begin{equation*}
E[\Lambda ^{-p}]\leq C_{5}+C_{6}\int_{0}^{1}\epsilon ^{-k^{\prime }}P\left(
N_{\epsilon }(t_{0})>\left\lfloor t_{0}\left( \frac{\epsilon _{0}}{\epsilon }%
\right) ^{1/K(q)}\right\rfloor \right) d\epsilon .
\end{equation*}%
From the proof of theorem \ref{cthm} we see that $K(q)=K(q,\epsilon )=\beta
^{-1}$ for $\epsilon $ small enough, where $\beta <m(j_{0})$, and hence to
see that $E[\Lambda ^{-p}]<\infty $ it will suffice to show 
\begin{equation*}
P\left( N_{\epsilon }(t_{0})>\left\lfloor t_{0}\left( \frac{\epsilon _{0}}{%
\epsilon }\right) ^{\beta }\right\rfloor \right) \text{ \ is \ }o(\epsilon
^{q})\text{ as }\epsilon \rightarrow 0\text{ for any }q>0.
\end{equation*}%
\noindent Chebyshev's inequality and (\ref{rategrowth}) yield 
\begin{eqnarray*}
P\left( N_{\epsilon }(t_{0})>\left\lfloor t_{0}\left( \frac{\epsilon _{0}}{%
\epsilon }\right) ^{\beta }\right\rfloor \right) &\leq &\exp \left(
-t_{0}\left( \frac{\epsilon _{0}}{\epsilon }\right) ^{\beta
}+(e-1)t_{0}\lambda (\epsilon )\right) \\
&\leq &\exp \left( -t_{0}\left( \frac{\epsilon _{0}}{\epsilon }\right)
^{\beta }+C(e-1)t_{0}f(\epsilon )\right) \text{ as }\epsilon \rightarrow 0.
\end{eqnarray*}

\noindent Which, by the definition of $f$ is seen to be $o(\epsilon ^{q})$
for any $q>0$ if%
\begin{equation*}
\beta >\frac{3\delta (\kappa -n)}{(\kappa -n+\alpha )}.
\end{equation*}%
Since $\beta $ and $\delta $ may take any values subject to the constraints $%
\beta <m(j_{0})$ and

\noindent $16\delta >\max \left( 8-r+v/2),(\kappa -n+\alpha )/4\alpha
\right) ,$ this condition becomes%
\begin{equation*}
16m(j_{0})>3(\kappa -n)\max \left( \frac{8-r+v/2}{\kappa -n+\alpha },\frac{1%
}{4\alpha }\right) .
\end{equation*}
\end{proof}

\bigskip

The condition (\ref{bracketcond}) exposes the qualitative structure of the
problem structure of the problem quite well in that it becomes easier to
satisfy with smaller values of $j_{0}$ (so that $%
%TCIMACRO{\U{211d} }%
%BeginExpansion
\mathbb{R}
%EndExpansion
^{e}$ is spanned with brackets of smaller length), or with smaller values of 
$\kappa $ (less intense jumps) or larger values of $\alpha $ (corresponding
to better behaved vector fields). One might think that the use of the lower
bound $t_{0}(m+1)^{-1}$ on the size of the longest interval is somewhat
crude. Indeed, conditional on $N_{\epsilon }(t_{0})=m,$ the distribution
function of the longest interval is known (see Feller \cite{feller}) :%
\begin{equation*}
F(x)=\overset{m}{\underset{i=1}{\sum }}(-1)^{-i}\dbinom{m}{i}\left( 1-\frac{%
ix}{t_{0}}\right) _{+}^{i-1}
\end{equation*}

\noindent and more explicit calculation may be performed using this, however
they seem to lead to no improvement in the eventual criterion obtained.
Clearly, the use of only part of the covariance matrix in forming the
estimate is an area in which improvement would allow further insight to be
gained.

\end{document}